# Linear parameterization of group GL(N, C)

© A. Lavrenov

E-mail: lanin99@mail.ru

**Abstract:** A linear parameterization of group GL(N, C) formed by direct products of matrices with in advance known symmetry properties is offered. Initial conditions of the given approach and the proposal on realization for $N = 2^m$ are discussed. The concrete example of application of the above-stated approach to group GL(4, C) is considered.

## 1 Introduction

The application of groups is useful for finding-out symmetry properties of studied systems. Therefore consideration in the given work of general linear group GL(N, C) is an actual problem. Input of effective parameterization for such group allows to investigate invariant characteristics of the composite law at group level. The complex vector parameterization of Lorentz group in [1] can serve as an evident and indicative example. Attempt to extend or enter similar to parameterization on others, great values of an order of group N were undertaken earlier [2-4]. However their efficiency and pragmatic are sufficient are small due to various reasons – bulkiness and a fog of calculations, immensity of qualitative results etc. Therefore convenient and evident parameterization for the above-stated group remains the open problem demanding of the decision presently. In the given work the perspective approach to linear parameterization of group GL(N, C) formed by certain basis of matrices with in advance known symmetry properties is offered. In the first part initial conditions of the given approach and the proposal on realization are discussed. In the second part the concrete example of application of the above-stated approach to group GL(4, C) is offered. The final part is devoted discussion of possible generalizations and directions of the further researches.

## 2 Preliminary data and the proposal on realization

It is known that any matrix in the size $N \times N$ can be spread out on basis direct (kroneker) multiplication of basic vectors of $N$-dimensional space:

$$A = \sum_{i,j=1}^{i,j=N} A_{ij} e_i \otimes e_j \equiv A_{ij} e_{ij} \tag{1}$$

Let's name such numbering for elements of a matrix as global, and representation of a matrix – as standard. However there are other bases and accordingly representations of matrices. Advantage for use of any concrete matrix basis is defined by symmetry of a considered problem. With that end in view we will remember different systems of coordinates (SC) on which it is possible to spread out, for example, a vector. Its coordinates generally have any values of components, and it, accordingly, possesses symmetry SC or not. A simple example of the given remark for matrices is input of symmetric $e_{ij}^s$ and antisymmetric $e_{ij}^a$ bases:

$$e_{ij}^s = \frac{e_{ij} + e_{ji}}{2} \quad \text{and} \quad e_{ij}^a = \frac{e_{ij} - e_{ji}}{2}$$

Then we receive for matrix **A** following decomposition:

$$A = A_{ij} e_{ij} = A_{ij} (\frac{e_{ij} + e_{ji}}{2} + \frac{e_{ij} - e_{ji}}{2}) = A_{ij} (e_{ij}^s + e_{ij}^a)$$

If matrix **A** has components only in one of bases $e_{ij}^s$ ($e_{ij}^a$), then it will be a symmetric (antisymmetric) matrix.

Thus, it is possible to formulate the following statement:

**THE STATEMENT 1.** *Each element of matrix basis defines the type of matrixes or their symmetry.*

Let's transform Eq. (1), allocating in an explicit form separate items:

$$A = \sum_{i=1}^{i=N} \sum_{j=1}^{j=N} A_{ij} e_i \otimes e_j \equiv \sum_{i=1}^{i=i'} \sum_{j=1}^{j=j'} A_{ij} e_i \otimes e_j + \sum_{i=i'+1}^{i=N} \sum_{j=1}^{j=j'} A_{ij} e_i \otimes e_j + \sum_{i=1}^{i=i'} \sum_{j=j'+1}^{j=N} A_{ij} e_i \otimes e_j + \sum_{i=i'+1}^{i=N} \sum_{j=j'+1}^{j=N} A_{ij} e_i \otimes e_j$$

Them we consider as independent objects where it is possible to enter other numbering. We will name its as local numbering. Such transformation means decomposition of an initial matrix on separate, certain set of the matrices, which represent a part of matrix basis. In this case local numbering is transfer for elements of an initial matrix in aforementioned set of matrices:

$$A = \sum_{i=1}^{i=N}\sum_{j=1}^{j=N} A_{ij} e_i \otimes e_j \equiv \sum_{k=1}^{k=i'}\sum_{l=1}^{l=j'} A_{kl} e_k \otimes e_l + \sum_{m=1}^{m=N-i'}\sum_{l=1}^{l=j'} A_{ml} e_m \otimes e_l + \sum_{k=1}^{k=i''}\sum_{n=1}^{n=N-j'} A_{kl} e_k \otimes e_n + \sum_{m=1}^{m=N-i'}\sum_{n=1}^{n=N-j'} A_{mnj} e_m \otimes e_n$$

or

$$A = A_{ij} e_i \otimes e_j \equiv A_{kl} e_{kl} + A_{ml} e_{ml} + A_{kn} e_{kn} + A_{mn} e_{mn}.$$

In last equality global numbering is carried out by indexes *i, j,* and local – indexes *k, l, m, n*. Communication between indexes is easily established as *i=k* (*j=l*) at *j<l'+1* (*j<j'+1*) and *i=i'+m, j=j'+n* otherwise.

For a qualitative example it is possible to consider a situation of closed subspaces which give block-diagonal structure of an initial matrix. In this case zero values of elements of an initial matrix correspond $A_{ml} = A_{kn} \equiv 0$ in local numbering.

Thus, it is possible to formulate two more statements:

**THE STATEMENT 2.** *Allocation in initial N-dimensional space (not) equal subspaces (in other words, certain basis vectors) conducts to allocation in a matrix of corresponding (rectangular) square blocks. Therefore, we will name such representation of an initial matrix as block.*

Let's make two remarks: 1) the example considered above is the elementary case for presence of two subspaces; 2) the inconsistent choice of basis vectors is template allocation for elements of an initial matrix which after corresponding redefinition of basis vectors, is reduced to the block.

**THE STATEMENT 3.** *Global and local numberings for elements of an initial matrix in block representation are connected by linear transformation.*

Let's discuss now properties of direct multiplication at the description of matrices. That fact in this case is meant that certain types of matrices or their structure can be described by means of direct product of initial any matrices.

Hence, then according to the statement 1 we receive incorporated type of matrices which has symmetry of initial matrices. This fact we fix in the following statement.

**THE STATEMENT 4.** *Possibility of record for an initial matrix in the form of direct multiplication for initial any matrices conducts to possibility to allocation in it of corresponding templates. The template form is defined by types of the matrices entering into direct product. We will name such representation for an initial matrix as template.*

Lexicographic numbering which allows to connect usual, global numbering of an initial matrix with local numbering in the matrices entering into direct product, gives the chance to see visually interrelation of allocation of certain templates in an initial matrix. In other words, allocation of certain subspaces of the dimensions $n^{(k)}$ connected unequivocally with the size of matrixes of direct product, in initial *N*-dimensional space is available. So if to designate

$$A = A_{ij} e_{ij}; A' = A'_{i'j'} e_{i'j'}; A'' = A''_{i''j''} e_{i''j''};\ldots; A^{(n)} = A_{i^{(n)} j^{(n)}} e_{i^{(n)} j^{(n)}},$$

then $A = A_{ij} e_{ij} \equiv (A^{(n)} \otimes \ldots \otimes A'' \otimes A')_{ij} e_{ij} \equiv A^{(n)}{}_{i^{(n)} j^{(n)}} \ldots A''_{i''j''} A'_{i'j'} e_{i^{(n)} j^{(n)}} \otimes e_{i''j''} \otimes e_{i'j'},$

where

$i = n^{(n-1)} \ldots n'' n'(i^{(n)} - 1) + \ldots + n'(i'' - 1) + i'$ and $j = n^{(n-1)} \ldots n'' n'(j^{(n)} - 1) + \ldots + n'(j'' - 1) + j'$.

**THE STATEMENT 5.** *Global and local numberings for elements of an initial matrix in the template representation caused by direct product are connected by nonlinear transformation.*

Completeness and uniqueness for the description of an initial matrix remain former in different representations thanks to the unequivocal conformity of local and global numbering, i.e. all above described representations of matrices are identical among themselves. Novelty, advantage of template representation appears at use of such matrices which will play in our approach, on the one hand, a role base, and on the other hand, will possess known symmetry properties with the linear law of a composition. In this case matrices of the big sizes can be parameterized effectively by consideration of a composition for base matrices of a smaller size. In the universal case for these purposes we take so-called the generalized Gell-Mann matrices with their known linear law of a composition. However now we consider the most elementary case of Pauli's matrices $\sigma_\mu = (\sigma_0, \vec{\sigma}) = (\sigma_0, \sigma_k) = (\sigma_0, \sigma_1, \sigma_2, \sigma_3)$ with their corresponding symmetry properties:

$$\sigma_i \sigma_j = \delta_{i,j} + i\varepsilon_{ijk}\sigma_k; \qquad \sigma_0\sigma_k = \sigma_k\sigma_0 = \sigma_k; \sigma_0\sigma_0 = I_2; \qquad Tr(\sigma_k) = 0, \qquad (2)$$

where $\sigma_0 = \begin{pmatrix} 1 & 0 \\ 0 & 1 \end{pmatrix}; \quad \sigma_1 = \begin{pmatrix} 0 & 1 \\ 1 & 0 \end{pmatrix}; \quad \sigma_2 = \begin{pmatrix} 0 & -i \\ i & 0 \end{pmatrix}; \quad \sigma_3 = \begin{pmatrix} 1 & 0 \\ 0 & -1 \end{pmatrix}; \quad I_n$ - $n$-dimensional unit matrix.

Let's remind that 1) by means of Pauli's matrices completely and uniquely it is possible to present any two-dimensional (2x2) matrix; 2) their use for rotation group already showed the efficiency [1].

Thus, our offer we formulate in the form of following theorems.

**THEOREM 1.** *Any square matrix in size $N \times N$ ( $N = 2^m$ because of a choice of Pauli's matrices) can be spread out uniquely in basis $\sigma_{\mu_1} \otimes \ldots \otimes \sigma_{\mu_m}$:*

$$A = A_{ij}e_{ij} \equiv A_{\mu_1\ldots\mu_m}\sigma_{\mu_1} \otimes \ldots \otimes \sigma_{\mu_m} \equiv A_m\sigma^{[m]} = A_{L\mu_k R}\sigma_L^{[k-1]} \otimes \sigma_{\mu_k} \otimes \sigma_R^{[m-k+1]} \qquad (3)$$

*where designations are entered*

$$\sigma^{[m]} = \sigma_{\mu_1} \otimes \ldots \otimes \sigma_{\mu_m}; \qquad \sigma_L^{[k-1]} = \sigma_{\mu_1} \otimes \ldots \otimes \sigma_{\mu_{k-1}}; \qquad \sigma_R^{[m-k+1]} = \sigma_{\mu_{k+1}} \otimes \ldots \otimes \sigma_{\mu_m}.$$

The theorem proof is easy for executing on a basis application of property for Pauli's matrices (2) and direct product:

$$Tr(A \otimes B) = Tr(A)Tr(B)$$

**THEOREM 2.** *The law of a composition for any matrices of the theorem 1 (i.e. presented to basis $\sigma_{\mu_1} \otimes \ldots \otimes \sigma_{\mu_m}$), is linear.*

To prove this theorem, we present for the composite law a following equation which allows to write out certain item **C=AB** for each concrete case. We allocate in decomposition (3), for example, member $\Pi_L \otimes \sigma_k \otimes \Pi_R$, where $\Pi_L, \Pi_R$ can be empty elements. Then the composition or multiplication of such members conducts according to (2) to following result

$$(\Pi_L^A \otimes \sigma_i \otimes \Pi_R^A)(\Pi_L^B \otimes \sigma_j \otimes \Pi_R^B) =$$
$$= \delta_{ij} \Pi_L^A \Pi_L^B \otimes I_2 \otimes \Pi_R^A \Pi_R^B + i\varepsilon_{ijk} \Pi_L^A \Pi_L^B \otimes \sigma_k \otimes \Pi_R^A \Pi_R^B \qquad (4)$$

## 3 Group GL(4, C)

According to our theorem in basis $\sigma_{\mu_1} \otimes \sigma_{\mu_2}$, any matrix **A** of the 4x4 sizes looks as follows

$$A = A_{\mu_1\mu_2}\sigma_{\mu_1} \otimes \sigma_{\mu_2} \equiv A_{00}I_2 \otimes I_2 + A_{i0}\sigma_i \otimes I_2 + A_{0j}I_2 \otimes \sigma_j + A_{kl}\sigma_k \otimes \sigma_l. \qquad (5)$$

Hence, the composite law will look as follows:

$$C = AB = C_{\mu_1\mu_2}\sigma_{\mu_1} \otimes \sigma_{\mu_2} \equiv C_{00}I_2 \otimes I_2 + C_{k0}\sigma_k \otimes I_2 + C_{0l}I_2 \otimes \sigma_l + C_{ij}\sigma_i \otimes \sigma_j, \qquad (6)$$

where

$$C_{00} = A_{00}B_{00} + A_{j0}B_{j0} + A_{0j}B_{0j} + A_{kl}B_{kl};$$

$$C_{k0} = A_{00}B_{k0} + A_{k0}B_{00} + A_{0j}B_{kj} + A_{kj}B_{0j} + i\varepsilon_{lmk}[A_{l0}B_{m0} + A_{lj}B_{mj}];$$

$$C_{0l} = A_{00}B_{0l} + A_{0l}B_{00} + A_{il}B_{i0} + A_{i0}B_{il} + i\varepsilon_{jkl}[A_{0j}B_{0k} + A_{ij}B_{ik}];$$

$$C_{ij} = A_{00}B_{ij} + A_{0j}B_{i0} + A_{i0}B_{0j} + A_{ij}B_{00} + i\varepsilon_{klj}[A_{ik}B_{0l} + A_{ok}B_{il}] + i\varepsilon_{kli}[A_{kj}B_{l0} + A_{k0}B_{lj}] - \varepsilon_{kmi}\varepsilon_{\ln j}A_{kl}B_{mn}$$

The given composite law allows to survey at once a class of such matrices which are closed concerning multiplication. In particular, the account of the first and second or third members in decomposition (5) matrix **A** makes evident such closed class. The consideration of the law for a composition of antisymmetric matrices **A** can serve as other, more difficult example – example for open class. For these purposes we will remind that

$$(A \otimes B)^T = A^T \otimes B^T; \qquad \sigma_0^T = \sigma_0; \qquad \sigma_1^T = \sigma_1; \qquad \sigma_2^T = -\sigma_2; \qquad \sigma_3^T = \sigma_3, \qquad (7)$$

where T – a transposing sign.

Therefore, the antisymmetric matrix **A** has the following structure taking into account (7) generally

$$A = A_{20}\sigma_2 \otimes I_2 + A_{21}\sigma_2 \otimes \sigma_1 + A_{23}\sigma_2 \otimes \sigma_3 + A_{02}I_2 \otimes \sigma_2 + A_{12}\sigma_1 \otimes \sigma_2 + A_{32}\sigma_3 \otimes \sigma_2, \qquad (8)$$

and the law for a composition of antisymmetric matrixes **A** and **B** will look so

$$C = AB = C_{00}I_2 \otimes I_2 + C_{i0}\sigma_i \otimes I_2 + C_{0j}I_2 \otimes \sigma_j + C_{kl}\sigma_k \otimes \sigma_l,$$

where

$C_{0,0} = A_{02}B_{02} + A_{12}B_{12} + A_{23}B_{23} + A_{20}B_{20} + A_{21}B_{21} + A_{32}B_{32};$ $\quad C_{0,1} = A_{20}B_{21} + A_{21}B_{20};$

$C_{0,2} = i(A_{23}B_{21} - A_{21}B_{23});$ $\quad C_{0,3} = A_{20}B_{23} + A_{23}B_{20};$ $\quad C_{1,0} = A_{02}B_{12} + A_{12}B_{02};$

$C_{2,0} = i(A_{32}B_{12} - A_{12}B_{32});$ $\quad C_{3,0} = A_{02}B_{32} + A_{32}B_{02};$ $\quad C_{1,1} = A_{23}B_{32} + A_{32}B_{23};$

$C_{1,2} = i(A_{20}B_{32} - A_{32}B_{20});$ $\quad C_{1,3} = -A_{21}B_{32} - A_{32}B_{21};$ $\quad C_{2,1} = i(A_{02}B_{23} - iA_{23}B_{02});$

$C_{2,2} = A_{02}B_{20} + A_{20}B_{02};$ $\quad C_{2,3} = i(A_{21}B_{02} - A_{02}B_{21});$ $\quad C_{3,1} = -A_{12}B_{23} - A_{23}B_{12};$

$C_{3,2} = i(A_{12}B_{20} - A_{20}B_{12});$ $\quad C_{3,3} = A_{12}B_{21} + A_{21}B_{12}.$

In [1] antisymmetric matrix **A** was parameterized by a complex vector $\vec{q} = \vec{a} + i\vec{b}$:

$$A = \begin{bmatrix} a^\times & i\vec{b} \\ -i\vec{b} & 0 \end{bmatrix}, \quad (9)$$

where $(a^\times)_{ij} = \varepsilon_{ijk}a_k$.

Comparing formulas (8) and (9), it is easy to find relation between components $A_{\mu_1\mu_2}$ and $q_i = a_i + ib_i$:

$-ia_1 = A_{21} - A_{12};$ $\quad -ia_2 = -A_{20} - A_{23};$ $\quad -ia_3 = A_{02} + A_{32};$

$b_1 = -A_{21} - A_{12};$ $\quad b_2 = -A_{20} + A_{23};$ $\quad b_3 = -A_{02} + A_{32}.$

## 4 Conclusion

In the given work the perspective approach to the linear parameterization of group GL(N, C) formed by direct products of matrices with in advance known symmetry properties is offered. Initial conditions of the given approach and the proposal on realization for $N = 2^m$ are discussed. The concrete example of application of the above-stated approach to group GL(4, C) is considered. The further directions of researches can be designated as specification of distinctions in bases and them symmetry properties for classification problems of matrices of any dimension, and also practical application of the offered approach in various problems.